\documentclass{article}
\usepackage{amsthm}
\usepackage{amsmath,a4wide}
\usepackage{amssymb,color}

\theoremstyle{theorem}
\newtheorem{theorem}{Theorem}

\theoremstyle{definition}

\newtheorem*{remark}{Remark}

\newtheorem{lemma}{Lemma}
\newtheorem{coro}{Corollary}

\def\ni{\noindent}

\begin{document}
\title{Power Partitions and Semi-$m$-Fibonacci Partitions}
\date{}
\maketitle
\begin{center}
Abdulaziz M. Alanazi$^{1}$, Augustine O. Munagi$^{2}$, Darlison Nyirenda$^{3}$\\
$^{1}$Department of Mathematics, Faculty of Sciences, University of Tabuk, P.O. Box 741, Tabuk 71491, Saudi Arabia\\
$^{2,3}$School of Mathematics, University of the Witwatersrand, P.O. Wits 2050, Johannesburg, South Africa \\
$^{1}$am.alenezi@ut.edu.sa, $^{2}$augustine.munagi@wits.ac.za, $^{3}$darlison.nyirenda@wits.ac.za
\end{center}
 
\begin{abstract}
George Andrews recently proved a new identity between the cardinalities of the set of Semi-Fibonacci partitions and the set of partitions into powers of two with all parts appearing an odd number of times. This paper extends the identity to the set of Semi-$m$-Fibonacci partitions of $n$ and the set of partitions of $n$ into powers of $m$ in which all parts appear with multiplicity not divisible by $m$. We also give a new characterization of Semi-$m$-Fibonacci partitions and some congruences satisfied by the associated number sequence. 
\bigskip

\smallskip

\end{abstract}

\section{Introduction}\label{introd}
A partition $\lambda$ of an integer $n>0$ is a finite nonincreasing integer sequence whose sum is $n$. The terms of the sequence are called \emph{parts} of $\lambda$.  
Thus a partition with $k$ parts will generally be expressed as
\begin{equation}\label{notn1}  
\lambda = (\lambda_1,\lambda_2,\dots,\lambda_k),\ \lambda_1\ge \lambda_2\ge \cdots\ge \lambda_k > 0,
\end{equation}
or
\begin{equation}\label{notn2}  
\lambda = (\lambda_1^{v_1},\lambda_2^{v_2},\dots,\lambda_t^{v_t}),\ \lambda_1 > \lambda_2 > \cdots > \lambda_t > 0,\, t\le k,
\end{equation}
where $\lambda_i^{v_i}$ indicates that $\lambda_i$ occurs with multiplicity $v_i$, for each $i$, and $v_1+\cdots + v_t = k$ \cite{A1}. 

In a recent paper paper Andrews \cite{A0} describes the set $SF(n)$ of semi-Fibonacci partitions as follows: $SF(1)=\{(1)\},\, SF(2)=\{(2)\}$. If $n>2$ and $n$ is even then

$SF(n)=\{\lambda\mid \lambda\ \mbox{is a partition of $\frac{n}{2}$ with each part doubled}\}$. If $n$ is odd, then a member of $SF(n)$ is obtained by inserting 1 into each partition in $SF(n-1)$ or by adding 2 to the single odd part in a partition in $SF(n-2)$.  

The cardinality $sf(n)=|SF(n)|$ satisfies the following recurrence relation for all $n>0$ (with $sf(-1)=0, sf(0)=1$);
\begin{equation}\label{eq00}
sf(n) = \begin{cases}
sf(n/2),& \mbox{if $n$ is even};\\
sf(n-1)+sf(n-2),& \mbox{if $n$ is odd}.
\end{cases}
\end{equation}

The semi-Fibonacci sequence $\{sf(n)\}_{n>0}$ occurs as sequence number A030067 in Sloane's database \cite{Sl}. George Beck \cite{B} has previously considered the properties of a set of polynomials related to the semi-Fibonacci partitions. 

Andrews stated the following relation between the number of semi-Fibonacci partitions of $n$ and the number $ob(n)$ of binary partitions of $n$ in which every part occurs an odd number of times:
\begin{theorem}[\cite{A0}, Theorem 1]\label{thm0} For each $n\geq 0$,
\begin{equation}\label{eq01}
sf(n) = ob(n),
\end{equation}
\end{theorem}
Andrews gave a generating function proof and asked for a bijective proof.

The proof turns out to be remarkable simple. It goes as follows.
Each part $t$ of $\lambda\in SF(n)$ can be expressed as $t=2^i\cdot h,\, i\geq 0$, where $h$ is odd. Now transform $t$ as
$$t=2^i\cdot h\longmapsto 2^i,2^i,\ldots,2^i (h\ \mbox{times}).$$
This gives a partition of $n$ into powers of 2 in which every part has odd multiplicity.
Conversely, consider $\beta\in OB(n)$. Since every part (a power of 2) has odd multiplicity we simply write $\beta$ in the exponent notation $\beta=(\beta_1^{u_1},\ldots,\beta_s^{u_s}),\, \beta_1>\cdots >\beta_s$ with the $u_i$ odd and positive. Since each $\beta_i^{u_i}$ has the form $(2^{j_i})^{u_i},\, j_i\geq 0$, we apply the transformation:
$$\beta_i^{u_i}=(2^{j_i})^{u_i}\longmapsto 2^{j_i}u_i.$$
This gives a unique partition in $SF(n)$. Indeed the image may contain at most one odd part which occurs precisely when $j_i=0$.
\smallskip

\begin{table}[htp]
\begin{center}
\begin{tabular}{ccc}
$SF(9)$ & $\stackrel{}{\longrightarrow}$& $OB(9)$\\ \hline
  (8,1) & $\mapsto$ & (8,1)\\
  (4,3,2) & $\mapsto$ & (4,2,1,1,1)\\
  (6,3) & $\mapsto$ & (2,2,2,1,1,1)\\
  (5,4) & $\mapsto$ & (4,1,1,1,1,1)\\
  (7,2) & $\mapsto$ & (2,1,1,1,1,1,1,1)\\
  (9) & $\mapsto$ & (1,1,1,1,1,1,1,1,1)\\ \hline
\end{tabular}
\caption{The map $SF(n)\rightarrow OB(n)$ for $n=9$.}\label{tab1}
\end{center}
\end{table}

We also consider the following congruence which Andrews proved with generating functions.

\begin{theorem}[\cite{A0}, Theorem 2]\label{thm01} For each $n\geq 0$, $sf(n)$ is even if $3|n$ and odd otherwise.
\end{theorem}
\begin{proof} We give a combinatorial proof based on mathematical induction. The result holds for $n=1,2,3$ since $sf(1)=1=sf(2)$ and $sf(3)=|\{(1,2),(3)\}|=2$.
Now let $n>3$ and assume that the result holds for all integers less than $n$.

If $n\equiv 1\pmod{3}$, then $sf(n)$ is the sum of $sf(n-1)$ and $sf(n-2)$ which have opposite parities since, by the inductive hypothesis, $sf(n-1)$ is even (since $3|(n-1)$) and $sf(n-2)$ is odd.\\
If $n\equiv 2\pmod{3}$, then $sf(n)$ is the sum of $sf(n-1)$ which is odd (since $3\nmid(n-1)$) and $sf(n-2)$ is even. Thus $sf(n)$ is odd.\\
If $3|n$ and $n$ is even, then $sf(n)=sf(n/2)$. Since $3|\frac{n}{2}$, it follows that $sf(n/2)$ is even by the inductive hypothesis.
Lastly, if $3|n$ and $n$ is odd, then $sf(n)=sf(n-1)+s(n-2)$ which is even since $3\nmid (n-1)$ and $3\nmid (n-2)$.

This completes the proof.
\end{proof}

The following result is easily deduced from the definition of sets counted by $sf(n)$.
\begin{coro}\label{coro1} Given a nonnegative integer $v$,
 $$sf(2^v)=1.$$
\end{coro}

In Section \ref{mfibo} we define the semi-$m$-Fibonacci partitions by extending the previous construction using a fixed integer modulus $m>1$. A generalized identity is then stated  between the set of semi-$m$-Fibonacci partitions and the set of partitions into powers of $m$ with multiplicities not divisible by $m$ (Theorem \ref{thm1}). Then in Subsection \ref{character} we give an independent characterization of the semi-$m$-Fibonacci partitions.
Lastly, in Section \ref{congruence} we discuss some arithmetic properties satisfied by the semi-$m$-Fibonacci sequence.

\section{Generalization}\label{mfibo}
We generalize the set of semi-Fibonacci Partitions to the set $SF(n,m)$ of semi-$m$-Fibonacci Partitions as follows:

$SF(n,m)=\{(n)\},\ n=1,2,\ldots,m$

If $n>m$ and $n$ is a multiple of $m$, then 

$SF(n,m)=\{\lambda\mid \lambda\ \mbox{is a partition of $\frac{n}{m}$ with each part multiplied by}\ m\}$. 

If $n$ is not a multiple of $m$, that is, $n\equiv r\pmod{m},\, 1\leq r\leq m-1$, then $SF(n,m)$ arises from two sources: first, partitions obtained by inserting $r$ into each partition in $SF(n-r,m)$, and second, partitions obtained by adding $m$ to the single part of each partition $\lambda \in SF(n-m,m)$ which is congruent to $r\pmod{m}$ (since $\lambda$ contains exactly one part which is congruent to $r$ modulo $m$, see Lemma \ref{lem001} below).

\begin{lemma}\label{lem001}
Let $\lambda\in SF(n,m)$. 

If $m\mid n$, then every part of $\lambda$ is a multiple of $m$. 

If $n\equiv r\pmod{m},\, 1\leq r <m$, then $\lambda$ contains exactly one part $\equiv r\pmod{m}$.
\end{lemma}
\begin{proof}
If $m\mid n$, the parts of a partition in $SF(n,n)$ are clearly divisible by $m$ by construction. 

For induction note that $SF(r,m)=\{(r)\},\, r=1,\ldots,m-1$, so the assertion holds trivially. Assume that the assertion holds for the partitions of all integers $<n$ and consider $\lambda\in SF(n,m)$ with $1\leq r<m$. Then $\lambda$ may be obtained by inserting $r$ into a partition $\alpha\in SF(n-r,m)$. Since $\alpha$ consists of multiples of m (as $m|(n-r)$), $\lambda$ contains exactly one part $\equiv r\pmod{m}$. Alternatively $\lambda$ is obtained by adding $m$ to the single part of a partition $\beta\in SF(n-m,m)$ which is $\equiv r\pmod{m}$. Indeed $\beta$ contains exactly one such part by the inductive hypothesis. Hence the assertion is proved.
\end{proof}

\ni As an illustration we have the following sets for small $n$ when $m=3$:

$SF(1,3)=\{(1)\}$,

$SF(2,3)=\{(2)\}$,

$SF(3,3)=\{(3)\}$,

$SF(4,3)=\{(4),(3,1)\}$,

$SF(5,3)=\{(5),(3,2)\}$,

$SF(6,3)=\{(6)\}$,

$SF(7,3)=\{(7),(4,3),(6,1)\}$,

$SF(8,3)=\{(8),(5,3),(6,2)\}$,

$SF(9,3)=\{(9)\}$,

$SF(10,3)=\{(10),(6,4),(7,3),(9,1)\}$.\\
Thus if we define $sf(n,m)=|SF(n,m)|$, we obtain that 

$sf(1,3)=sf(2,3)=sf(3,3)=1,\ sf(4,3)=2,\ sf(5,3)=2,\ sf(6,3)=1,$ 

$sf(7,3)=3,\ sf(8,3)=3,\ sf(9,3)=1,\ sf(10,3)=4$.\\
Therefore, for $m>1$, we see that $sf(n,m)=0$ if $n<0$ and $sf(0,m)=1$, and for $n>0$,

\begin{equation}\label{eq1}
sf(n,m) = \begin{cases}
sf(n/m,m),& \mbox{if}\ n\equiv 0\pmod{m}; \\
sf(n-r,m)+sf(n-m,m),& \mbox{if}\ n\equiv r\pmod{m},0<r<m.
\end{cases}
\end{equation}
The case $m=2$ gives the function considered by Andrews: $sf(n,2)=sf(n)$.

Power partitions are partitions into powers of a positive integer $m$, also known as $m$-power partitions \cite{G}. Let $nd(n,m)$ be the number of $m$-power partitions of $n$ in which the multiplicity of each part is not divisible by $m$.
Thus, for example, $nd(10,3)=4$, the enumerated partitions being $(9,1),(3,3,1,1,1,1),(3,1,1,1,1,1,1,1),(1,1,1,1,1,1,1,1,1,1)$.

\begin{theorem}\label{thm1} For integers $n\geq 0,m>1$,
\begin{equation}\label{eq2}
sf(n,m) = nd(n,m),
\end{equation}
\end{theorem}
\begin{proof} We give two proofs, one analytic one combinatorial.
\smallskip 

\ni {First Proof}.
Let $|q| < 1$ and define 
\begin{equation}\label{eqproof0}
G_m(q) = \sum\limits_{n \geq 0}sf(n,m)q^{n}.
\end{equation}
Then we have
\begin{align}
G_m(q) & = \sum\limits_{n \geq 0}sf(mn,m)q^{mn} + \sum\limits_{n \geq 0}sf(mn + 1,m)q^{mn + 1} + \ldots + \sum\limits_{n \geq 0}sf(mn + m - 1,m)q^{mn + m - 1}\nonumber \\
    & = \sum\limits_{n \geq 0}sf(mn,m)q^{mn} + \sum\limits_{r = 1}^{m - 1}\sum\limits_{n \geq 0}sf(mn + r,m)q^{mn + r}\label{eqproof1}\\
    & = \sum\limits_{n \geq 0}sf(n,m)q^{mn} + \sum\limits_{r = 1}^{m - 1}\sum\limits_{n \geq 0}(sf(mn,m) + sf(mn + r - m,m))q^{mn + r} \,\,\,(\text{by}\,\,\eqref{eq1}) \nonumber\\
    & = \sum\limits_{n \geq 0}sf(n,m)q^{mn} + \sum\limits_{r = 1}^{m - 1}\sum\limits_{n \geq 0}(sf(n,m)q^{mn + r} + \sum\limits_{r = 1}^{m - 1}\sum\limits_{n \geq 0} sf(mn + r - m,m)q^{mn + r}.\nonumber \\
    & = (1 + \sum\limits_{r = 1}^{m - 1}q^{r})\sum\limits_{n \geq 0}sf(n,m)q^{mn}  + \sum\limits_{r = 1}^{m - 1}\sum\limits_{n \geq 0} sf(m(n - 1) + r,m)q^{mn + r}\nonumber\\
    & = G_m(q^{m})\sum\limits_{r = 0}^{m-1}q^{r} + \sum\limits_{r = 1}^{m - 1}\sum\limits_{n \geq 0} sf(mn + r,m)q^{mn + m + r} \nonumber\\
    & = G_m(q^{m})\sum\limits_{r = 0}^{m-1}q^{r} + q^{m}\sum\limits_{r = 1}^{m - 1}\sum\limits_{n \geq 0} sf(mn + r,m)q^{mn + r}\nonumber \\
    & = G_m(q^{m})\sum\limits_{r = 0}^{m-1}q^{r}) + q^{m}\left( \sum\limits_{n \geq 0}sf(n,m)q^{n} - \sum\limits_{n \geq 0} sf(mn,m)q^{mn} \right)\,\,\,(\text{by}\,\,\eqref{eqproof1})\nonumber \\
    & =G_m(q^{m})\sum\limits_{r = 0}^{m-1}q^{r} + q^{m}(G_m(q) - G_m(q^{m})) \nonumber\\
    & = \left( -q^{m} + \sum\limits_{r = 0}^{m - 1}q^{r}\right)G_m(q^{m}) + q^{m}G_m(q).
    \end{align}
Hence,
\begin{equation}\label{eq11}
G_m(q) = \frac{1 + q + q^{2} + q^{3} + \ldots + q^{m  - 1} - q^{m}}{1 - q^{m}}G_m(q^{m}).
\end{equation}
Equation \eqref{eq11} implies that
$$ G_m(q) =  \left(\frac{1 + q + q^{2}+ \ldots + q^{m  - 1} - q^{m}}{1 - q^{m}}\right)\left(\frac{1 + q + q^{2m} + \ldots + q^{(m  - 1)m} - q^{m^{2}}}{1 - q^{m^{2}}}\right)G_m(q^{m^{2}}),$$ and continuing the iteration, we get
$$ G_m(q) = \prod\limits_{n = 0}^{N}\left(\frac{1 + q^{m^{n}} + q^{2m^{n}} + \ldots + q^{(m - 1)m^{n}} - q^{m^{n + 1}}}{1 - q^{m^{n + 1}}}\right)G_m(q^{m^{N + 1}}).$$
Taking the limit as $N\rightarrow \infty$, we have $G_m(q^{m^{N + 1}}) \rightarrow G_m(0) = 1$ (since $|q| < 1$) so that
\begin{align*}
G_m(q)  & = \prod\limits_{n = 0}^{\infty}\left(\frac{1 + q^{m^{n}} + q^{2m^{n}} + \ldots + q^{(m - 1)m^{n}} - q^{m^{n + 1}}}{1 - q^{m^{n + 1}}}\right) \\
      & = \prod\limits_{n = 0}^{\infty}\left(\frac{q^{m^{n}} + q^{2m^{n}} + \ldots + q^{(m - 1)m^{n}} + 1 - q^{m^{n + 1}}}{1 - q^{m^{n + 1}}}\right) \\
      & = \prod\limits_{n = 0}^{\infty}\left(1 + \frac{q^{m^{n}} + q^{2m^{n}} + \ldots + q^{(m - 1)m^{n}} }{1 - q^{m^{n + 1}}}\right) \\
      & = \prod\limits_{n = 0}^{\infty}\left(1 + (q^{m^{n}} + q^{2m^{n}} + \ldots + q^{(m - 1)m^{n}})\sum\limits_{j = 0}^{\infty}q^{j(m^{n + 1})}\right).
  \end{align*}
Thus,
\begin{align}
G_m(q)  & = \prod\limits_{n = 0}^{\infty}\left( 1 + \sum\limits_{j = 0}^{\infty}q^{m^{n}(jm + 1)} + \sum\limits_{j = 0}^{\infty}q^{m^{n}(jm + 2)} + \sum\limits_{j = 0}^{\infty}q^{m^{n}(jm + 3)} + \ldots + \sum\limits_{j = 0}^{\infty}q^{m^{n}(jm + m - 1)} \right)\nonumber \\
      & = \sum\limits_{n \geq 0}nd(n,m)q^{n}.\label{eqproof2}
\end{align}
The assertion follows by comparing coefficients in \eqref{eqproof0} and \eqref{eqproof2}.
\medskip

\ni {Second Proof}.
Each part $t$ of $\lambda\in SF(n,m)$ can be expressed as $t=m^i\cdot h,\, i\geq 0$, where $m$ does not divide $h$. Now transform $t$ as
$$t=m^i\cdot h\longmapsto m^i,m^i,\ldots,m^i (h\ \mbox{times}).$$
This gives a partition of $n$ into powers of $m$ in which every part has multiplicity not divisible by $m$.
Conversely, consider $\beta\in ND(n,m)$. Since every part (a power of $m$) has a non-multiple of $m$ as multiplicity we simply write $\beta$ in the exponent notation $\beta=(\beta_1^{u_1},\ldots,\beta_s^{u_s}),\, \beta_1>\cdots >\beta_s$ with the $u_i\not\equiv 0\pmod{m}$. Since each $\beta_i^{u_i}$ has the form $(m^{j_i})^{u_i}$, we apply the transformation:
$$\beta_i^{u_i}=(m^{j_i})^{u_i}\longmapsto m^{j_i}u_i.$$
This gives a unique partition in $SF(n,m)$. If $m\mid n$, this image contains only multiples of $m$. If $n\equiv r\pmod{m},\, 1\leq r<m$, the image consists of multiples of $m$ and exactly one part $\equiv r$ (mod $m$) which occurs when $j_i=0$.
\end{proof}

\begin{table}[htp]
\begin{center}
\begin{tabular}{ccc}
$SF(11,3)$ & $\stackrel{}{\longrightarrow}$& $ND(11,3)$\\ \hline
  (11) & $\mapsto$ & (1,1,1,1,1,1,1,1,1,1,1)\\
  (8,3) & $\mapsto$ & (3,1,1,1,1,1,1,1,1)\\
  (6,5) & $\mapsto$ & (3,3,1,1,1,1,1)\\
  (9,2) & $\mapsto$ & (9,1,1)\\ \hline
\end{tabular}
\caption{The map $SF(n,m)\rightarrow ND(n,m)$ for $n=11,\, m=3$.}\label{tab2}
\end{center}
\end{table}

\subsection{A characterization of Semi-$m$-Fibonacci Partitions}\label{character}
Define the \emph{max $m$-power} of an integer $N$ as the largest power of $m$ that divides $N$ (not just the exponent of the power). Thus using the notation $x_m(N)$, we find that $N=u\cdot m^s,\, s\geq 0$, where $m\nmid u$ and $x_m(N)=m^s$. So $x_m(N)>0$ for all $N$. 

For example, $x_2(50)=2,\, x_2(40)=8,\, x_3(216)=27$ and $x_5(216)=1$.

Note that if the parts of a partition $\lambda$ have distinct max $m$-powers, then the parts are distinct. For if $u\cdot m^s =\lambda_i=\lambda_j = v\cdot m^t\in \lambda$ with $m\nmid u,v$, and $s>t$, then $u\cdot m^{s-t}=v\implies m|v$ a contradiction.  

We define three (reversible) operations on a partition  $\lambda=(\lambda_1,\ldots,\lambda_k)$ with an integer $m>1$:

(i) If the last part of $\lambda$ is less than $m$, delete it: $\tau_1(\lambda)=(\lambda_1,\ldots,\lambda_{k-1})$;

(ii) If $m\nmid \lambda_t>m$, then $\tau_2(\lambda)=
(\lambda_1,\ldots,\lambda_{t-1},\lambda_{t}-m,\lambda_{t+1},\ldots,\lambda_k)$.

(iii) If $\lambda$ consists of multiples of $m$, divide every part by $m$: $\tau_3(\lambda)=(\lambda_1/m,\ldots,\lambda_k/m)$.\\
These operations are consistent with the recursive construction of the set $SF(n,m)$, where $\tau_3^{-1}, \tau_1^{-1}$ and $\tau_2^{-1}$ correspond, respectively, to the three quantities in the recurrence \eqref{eq1}.

\begin{lemma}\label{lem2} Let $B(n,m)$ denote the set of partitions of $n$ in which the parts have distinct max $m$-powers and at most one non-multiple of $m$.
Then if $\lambda\in B(n,m)$ and $\tau_i(\lambda)\neq\emptyset$, then $\tau_i(\lambda)\in B(N,m),\, i=1,2,3$, for some $N$.
\end{lemma}
\begin{proof} Let $\lambda=(\lambda_1,\ldots,\lambda_k)\in B(n,m)$. If $\lambda$ contains one part less than $m$, the part is $\lambda_k$. So $\tau_1(\lambda)\in B(n-\lambda_k,m)$ since the max $m$-powers remain distinct. It is obvious that the parity of $\lambda$ is inherited by $\tau_2(\lambda)=
(\lambda_1,\ldots,\lambda_{t-1},\lambda_{t}-m,\lambda_{t+1},\ldots,\lambda_k)\in B(n-m,m)$.
Lastly, since the parts of $\lambda$ have distinct max $m$-powers $\tau_3(\lambda)=(\lambda_1/m,\ldots,\lambda_k/m)$ may contain at most one non-multiple of $m$ as a part. Hence $\tau_3(\lambda)\in B(n/m,m)$.
\end{proof}

We state an independent characterization of the Semi-$m$-Fibonacci Partitions.

\begin{theorem}\label{maxm}
A partition of $n$ is a semi-$m$-Fibonacci partition if and only if the parts have distinct max $m$-powers and at most one non-multiple of $m$. 
\end{theorem}

\begin{proof} We show that $SF(n,m)= B(n,m)$. Let $\lambda=(\lambda_1,\ldots,\lambda_k)\in SF(n,m)$ such that $\lambda\notin B(n,m)$. Assume that there are $\lambda_i > \lambda_j$ satisfying $x_m(\lambda_i)=x_m(\lambda_j)$ and let $\lambda_i=u_im^s,\, \lambda_j=u_jm^s$ with $m\nmid u_i,u_j$. Observe that $\tau_1$ deletes a part less than $m$ if it exists. So we can use repeated applications of $\tau_2$ to reduce a non-multiple modulo $m$, followed by $\tau_1$. This is tantamount to simply deleting the non-multiple of $m$, say $\lambda_t$, to obtain a member of $B(n-\lambda_t,m)$ from Lemma \ref{lem2}.
By thus successively deleting non-multiples, and applying $\tau_3^c,\, c>0$, we  
obtain a partition $\beta = (\beta_1,\beta_2,\ldots)$ with $\beta_i=v_im^w>\beta_j=v_jm^w$, where $m\nmid v_i,v_j$ and $w\leq s$. Then apply $\tau_3^w$ to obtain a partition $\gamma$ with two non-multiples of $m$. Then by Lemma \ref{lem001}, $\gamma\notin SF(n,m)$.
Therefore $\lambda\in SF(n,m)\implies \lambda\in B(n,m)$.

Conversely let $\lambda=(\lambda_1,\ldots,\lambda_k)\in B(n,m)$. If $\lambda=(t),\, 1\leq t\leq m$, then $\lambda\in SF(t,m)$. If $m|\lambda_i$ for all $i$, then $\tau_3(\lambda)= (\lambda_1/m,\ldots,\lambda_k/m)\in B(n/m,m)$ contains at most one part $\not\equiv 0\pmod{m}$, so $\lambda\in SF(n,m)$. Lastly assume that $n\equiv r\not\equiv 0\pmod{m}$. Then $r\in \lambda$ or $\lambda_t\equiv r\pmod{m}$ for exactly one index $t$. Thus $\tau_1(\lambda)=(\lambda_1,\ldots,\lambda_{k-1})$ consists of multiples of $m$ while $\tau_2(\lambda)=(\lambda_1,\ldots,\lambda_{t-1},\lambda_t-m,\lambda_{t+1},\ldots,\lambda_k)$ still contains one part $\not\equiv 0\pmod{m}$. In either case $\lambda\in SF(n,m)$. Hence $B(n,m)\subseteq SF(n,m)$. The the two sets are identical.
\end{proof}

\begin{remark}
Notice that Theorem \ref{maxm} certifies the second (bijective) proof of Theorem \ref{thm1}. If $\lambda=(\lambda_1,\ldots,\lambda_k)\in SF(n,m)$ but $\lambda\notin B(n,m)$ on account of having two parts $\lambda_i,\lambda_j$ such that $\lambda_i=u_im^s>\lambda_j=u_jm^s$ with $m\nmid u_i,u_j$, then it cannot have an inverse image. Assume that $\lambda$ maps to $\beta\in ND(n,m)$ which then includes the parts $m^{u_i+u_j}$ ($u_i+u_j$ copies of $m$). Then $u_i+u_j$ may be a multiple of $m$ (for example, when $u_i=1, u_j=m-1$) which implies that $\beta\notin ND(n,m)$, a contradiction. Alternatively the pre-image of $\beta$ would include the part $m(u_i+u_j)$ and so cannot be $\lambda$.  
\end{remark}

\section{Arithmetic Properties}\label{congruence}
We prove several congruence properties of the numbers $sf(n,m)$.

\begin{theorem}\label{conjsolv1} Let $n,m$ be integers with $n\geq 0,\, m>1$. Then  
$$sf(nm+1,m)=sf(nm+2,m)=\cdots =sf(nm+m-1,m) = \sum\limits_{j = 0}^{n}sf(j,m).$$
\end{theorem}
\begin{proof}
\noindent Let $J_{r,m}(q) = \sum\limits_{n\geq 0}sf(nm + r,m)q^{n}$ where  $ r = 1, 2, 3, \ldots m - 1$.
Then
\begin{align*}
J_{r,m}(q) & = \sum\limits_{n\geq 0}sf(nm,m)q^{n}  +  \sum\limits_{n\geq 0}sf(mn + r - m,m)q^{n} \,\,\,(\text{by} \,\,\eqref{eq1}) \\
        & =  \sum\limits_{n\geq 0}sf(n,m)q^{n}   + \sum\limits_{n\geq 0}sf(mn + r ,m)q^{n + 1}\\
        & = G_m(q) + q\sum\limits_{n\geq 0}sf(mn + r ,m)q^{n} \\
        & = G_m(q) + qJ_{r,m}(q)
\end{align*}
so that
\begin{equation}\label{conjsol}
J_{r,m}(q) = \frac{G_m(q)}{1 - q}.
\end{equation}
Since the right hand side of \eqref{conjsol} is independent of $r$, we must have $ J_{1,m}(q) = J_{2,m}(q) = \ldots = J_{m - 1,m}(q)$
so that $sf(nm+1,m)=sf(nm+2,m)=\cdots =sf(nm+m-1,m)$. Furthermore, from \eqref{conjsol}, we observe that
\begin{align*}
\sum\limits_{n \geq 0}sf(mn + r,m)q^{n} & = \sum\limits_{n\geq  0}q^{n} \sum\limits_{n \geq 0}sf(n,m)q^{n} \\
                                        & = \sum\limits_{n \geq 0} \sum\limits_{j = 0}^{n}sf(j,m)q^{n}
\end{align*}
which implies that $sf(mn + r,m) = \sum\limits_{j = 0}^{n}sf(j,m)$.
\end{proof}

\begin{coro}\label{corogen} Given integers $m\geq 2$, then for any $j\geq 0$ and a fixed $v\in\{0,1,\ldots,m\}$,
$$sf(m^j(mv+r),m)=v+1,\ 1\leq r\leq m-1.$$
\end{coro}
\begin{proof} By applying \eqref{eq1} several times (the case when $m\mid n$), it is clear that for any $j \geq 0$,  $sf(m^j(mv+r),m) = sf(m^{j - 1}(mv+r),m) = sf(m^{j - 2}(mv+r),m) = \ldots  = sf(mv + r,m)$. By the last equality in Theorem \ref{conjsolv1}, we have
$$ sf(mv + r, m)  = \sum\limits_{i = 0}^{v}sf(i,m) = 1 + \sum\limits_{i = 1}^{v}sf(i,m),\, v\geq 0,\, 1\leq r <m.$$
If  $1\leq v< m$, then $\sum\limits_{i = 1}^{v}sf(i,m) = \sum\limits_{i = 1}^{v}(sf(i - i,m) + sf(i - m, m))\ $ (by \eqref{eq1}). 
Since $0 < i \leq v < m$, we have $sf(mv + r, m) = 1 + \sum\limits_{i = 1}^{v} (1 + 0)  = 1 + v$. \\
If $v = m$, then $\sum\limits_{i = 1}^{v}sf(i,m) = \sum\limits_{i = 1}^{m - 1}sf(i,m) + sf(m,m) = m - 1 + sf(1,m) = m - 1 + 1 = m$; thus $sf(mv + r) = v + 1$ is true in this case. Finally, if $v = 0$, it is not difficult to see that $sf(r,m) = 1$.
\end{proof}

We note a few interesting special cases of Corollary \ref{corogen} below.
\begin{coro}\label{coro2} We have the following for any integer $m\geq 2$:

(i) $sf(m^i,m)=1,\ i\geq 0.$

(ii) $sf(m^ih,m)=1,\ 1\leq h\leq m-1,\, i\geq 0.$

(iii) Given an integer $n\geq 0$, then for each $n\in\{0,1,\ldots,m\}$,
$$sf(nm+1,m)=sf(nm+2,m)=\cdots =sf((n+1)m-1,m)=v+1.$$
\end{coro}
\begin{proof}
Part (i) is the case $h=1$ of part (ii). Parts (ii) and (iii) are obtained by setting $v=0$ and $j=0$,respectively, in Corollary \ref{corogen}.
\end{proof}
Note that part (i) of Corollary \ref{coro2} implies Corollary \ref{coro1}. Also when $m=2$, part (iii) gives just the three values $sf(1)=1$, $sf(3)=2$ and $sf(5)=3$, the parities of which are consistent with Theorem \ref{thm01}.
Part (iii) is a stronger version of Theorem \ref{conjsolv1} since the restriction of $n$ to the set $\{0,1,\ldots,m\}$ specifies a common value.

\begin{theorem}\label{thm6}
For any $j\geq 0$, $$\sum\limits_{r = 0}^{2j + 1}sf(r, 3) \equiv 0 \pmod{2}.$$
Consequently, 
\begin{equation}\label{cong61}
 sf(3j + 4,3) = sf(3j + 5,3) \equiv 0 \pmod{2}\,\,\text{where}\,\,j \equiv 0 \pmod{2}, 
\end{equation}
\begin{equation}\label{cong62}
sf(3^{r}j + 4,3) = sf(3^{r}j + 5,3) \equiv 0 \pmod{2}\,\, \text{for all}\,\, j\geq 0, r \geq 2.
\end{equation}
\end{theorem}
\begin{proof}
Note the following identity 
\begin{equation}\label{ident6}
\frac{1}{1 - q} = \prod\limits_{n = 0}^{\infty}(1 + q^{3^{n}} + q^{2\cdot 3^{n}}).
\end{equation}
Recall that
\begin{align*}
\sum\limits_{n \geq 0} sf(n,3)q^{n} & = \prod\limits_{n = 0}^{\infty}\left(\frac{1 + q^{3^{n}} + q^{2\cdot 3^{n}} - q^{3\cdot 3^{n}}}{1 -  q^{3\cdot 3^{n}}  }\right) \\
     & \equiv \prod\limits_{n = 0}^{\infty}\left(\frac{1 + q^{3^{n}} + q^{2\cdot 3^{n}} + q^{3\cdot 3^{n}}}{1 +  q^{3\cdot 3^{n}} }\right) \pmod{2}\\
     & = \prod\limits_{n = 0}^{\infty}\frac{(1 + q^{3^{n}})(1 + q^{2\cdot 3^{n}})}{1 +  q^{3\cdot 3^{n}} }\\
     & = \prod\limits_{n = 0}^{\infty}\left( \frac{ 1 + q^{2\cdot 3^{n}}}{1 +  q^{3^{n}} + q^{2\cdot 3^{n}} }\right) \\
     & = (1 - q)\prod\limits_{n = 0}^{\infty}(1 + q^{2\cdot 3^{n}}) \,\,\,(\text{by}\,\,\,\eqref{ident6}).     
\end{align*}
Thus
$$\frac{1}{1 - q}\sum\limits_{n \geq 0}sf(n,3)q^{n} \equiv \prod\limits_{n = 0}^{\infty}(1 + q^{2\cdot 3^{n}}) \pmod{2},$$
i.e.
$$\sum\limits_{n \geq 0}\sum\limits_{r = 0}^{n}sf(r,3)q^{n} \equiv \prod\limits_{n = 0}^{\infty}(1 + q^{2\cdot 3^{n}}) \pmod{2}.$$
Since the series expansion  of the right-hand side of the preceeding equation has even exponents, the result follows.\par \noindent
To prove \eqref{cong61}, we have
\begin{align*}
sf(3j + 4,3) & = sf(3(j + 1) + 1,3)\\
          &  = sf(3(j + 1) + 2,3) \,\,\,\,(\text{by Theorem}\,\,\ref{conjsolv1}) \\
          &  = \sum\limits_{r = 0}^{j + 1}sf(r,3)\,\,\,\,(\text{by Theorem}\,\,\ref{conjsolv1}) \\
          & \equiv 0 \pmod{2}\,\,\, (\text{since}\,\, j + 1 \,\,\text{is odd}).
\end{align*}
Furthermore, for \eqref{cong62}, observe that
$$ 3^{r-1}j + 1 \equiv 
\begin{cases}
0, & \text{if}\,\, j \equiv 1 \pmod{2}; \\
 1, & \text{otherwise}.
\end{cases}
 $$
Now, if $j$ is odd, then
\begin{align*}
sf(3^{r}j + 4,3) & = sf(3(3^{r - 1}j + 1) + 1,3) \\
              & = sf(3(3^{r -1}j + 1) + 2,3)\\
            & = \sum_{r = 0}^{3^{r -1}j + 1}sf(r,3) \,\,\,(\text{by Theorem}\,\,\ref{conjsolv1}) \\
            & = sf(3^{r - 1}j + 1,3) + \sum_{r = 0}^{3^{r -1}j}sf(r,3) \\
            & \equiv sf(3^{r -1}j + 1,3) \pmod{2}\,\,(\text{since}\,\, 3^{r -1}j \,\,\text{is odd})\\
             & = \sum_{r = 0}^{3^{r -2}j}sf(r,3)\\
            & \equiv 0 \pmod{2}\,\,\,(\text{since}\,\,3^{r -2}j\,\,\text{is odd}).
\end{align*}
\end{proof}

On the other hand, if $j$ is even, use \eqref{cong61}.

\begin{theorem}\label{thm7}
Let $k\equiv m+r\pmod{2m}$ and $k\le{m^2+r}$ for $1\leq r\leq m-1$. If $n\ge0$, $m\ge2$ and $n=m^ik$ for $i\ge0$, then $sf(n,m)$ is even. 
\end{theorem}

\begin{proof}
$k\equiv m+r\pmod{2m}$ and $k\le{m^2+r}$ for $1\leq r\leq m-1$ imply that $k=m(2t+1)+r\le{m^2+r}\Rightarrow 2t+1\le{m}$, for some positive integer $t$. Then from Corollary \ref{corogen}, we have
\begin{align*}
sf(m^ik,m) & = sf(m^i(m(2t+1)+r),m)  \\
        & =  sf(m(2t+1)+r,m) \,\,\,(\text{by} \,\,\eqref{eq1}) \\
        & = 2t+1+1 \,\,\,(\text{by Corollary \ref{corogen} and since} \,\, 2t+1\le m) \\
        & = 2t+2.
\end{align*}
\begin{remark}
When $m = 3$, Theorem \ref{thm7} reduces to Theorem \ref{thm6} without the restriction $k\le{m^2+r}$. 
\end{remark}

\end{proof}
\bibliographystyle{amsplain}

\end{document}